\documentclass[12pt,a4paper,english]{article}
\usepackage{amsmath,amssymb,mathrsfs}
\usepackage[T1]{fontenc}
\usepackage{float} 
\usepackage{verbatim}
\usepackage{booktabs}
\usepackage{graphicx}
\usepackage[margin=1.185in]{geometry}
\usepackage{bbold}
\usepackage{natbib}
\usepackage{babel}
\usepackage{setspace}\doublespace
\usepackage{subfigure}
\usepackage[labelsep=period,labelfont=bf]{caption}

\providecommand{\U}[1]{\protect\rule{.1in}{.1in}}

\def \tx{\textrm}

\usepackage{epstopdf}
\DeclareGraphicsExtensions{.pdf,.png,.jpg}

\begin{document}

\title{Bayesian Inference Methods for Univariate and Multivariate GARCH Models: a Survey}

\author{
  Audron\.{e} Virbickait\.{e}\\ 
\small{Universidad Carlos III de Madrid}\\ \small{ Getafe (Madrid), Spain, 28903}\\ \small{audrone.virbickaite@uc3m.es}  
 \\ \and M.\ Concepci\'on Aus\'in\\\small{Universidad Carlos III de Madrid}\\ \small{Getafe (Madrid), Spain, 28903}\\\small{ concepcion.ausin@uc3m.es} 
  \and Pedro Galeano\\\small{Universidad Carlos III de Madrid}\\ \small{Getafe (Madrid), Spain, 28903}\\ \small{pedro.galeano@uc3m.es}  
}

\date{}

\maketitle
\thispagestyle{empty}
\newpage
\begin{abstract}

This survey reviews the existing literature on the most relevant Bayesian inference methods for univariate and multivariate GARCH models. The advantages and drawbacks of each procedure are outlined as well as the advantages of the Bayesian approach versus classical procedures. The paper makes emphasis on recent Bayesian non-parametric approaches for GARCH models that avoid imposing arbitrary parametric distributional assumptions. These novel
approaches implicitly assume infinite mixture of Gaussian distributions on the standardized returns which have been shown to be more flexible and describe better the uncertainty about future volatilities. Finally, the survey
presents an illustration using real data to show the flexibility and usefulness of the non-parametric approach.

\end{abstract}

\thispagestyle{empty}
\newpage

\section{Introduction}

Understanding, modeling and predicting the volatility of financial time series has been extensively researched for more than 30 years and the interest in the subject is far from decreasing. Volatility prediction has a very wide
range of applications in finance, for example, in portfolio optimization, risk management, asset allocation, asset pricing, etc. The two most popular approaches to model volatility are based on the Autoregressive Conditional
Heteroscedasticity (ARCH) type and Stochastic Volatility (SV) type models. The seminal paper of \cite{Engle1982} proposed the primary ARCH model while \cite{Bollerslev1986} generalized the purely autoregressive ARCH into an ARMA-type model, called the Generalized Autoregressive Conditional Heteroscedasticity (GARCH) model. Since then, there has been a very large amount of research on the topic, stretching to various model extensions and
generalizations. Meanwhile, the researchers have been addressing two important topics: looking for the best specification for the errors and selecting the most efficient approach for inference and prediction.

Besides selecting the best model for the data, distributional assumptions for the returns are equally important. It is well known, that every prediction, in order to be useful, has to come with a certain precision measurement. In this way the agent can know the risk she is facing, i.e.\ uncertainty. Distributional assumptions permit to quantify this uncertainty about the future. Traditionally, the errors have been assumed to be Gaussian, however, it has been widely acknowledged that financial returns display fat tails and are not conditionally Gaussian. Therefore, it is common to assume a Student-t distribution, see \cite{Bollerslev1987}, \cite{He1999} and \cite{Bai2003}, among others. However, the assumption of Gaussian or Student-t distributions is rather restrictive. An alternative approach is to use a mixture of distributions, which can approximate arbitrarily any distribution given a sufficient number of mixture components. A mixture of two Normals was used by \cite{Bai2003}, \cite{Ausin2007} and \cite{Giannikis2008}, among others. These authors have shown that the models with the mixture distribution for the errors outperformed the Gaussian one and do not require additional restrictions on the degrees of freedom parameter as the Student-t one.

As for the inference and prediction, the Bayesian approach is especially well-suited for GARCH models and provides some advantages compared to classical estimation techniques, as outlined by \cite{Ardia2010}. Firstly,  the
positivity constraints on the parameters to ensure positive variance, may encumber some optimization procedures. In the Bayesian setting, constraints on the model parameters can be incorporated via priors. Secondly, in most of the cases we are more interested not in the model parameters directly, but in some non-linear functions of them. In the maximum likelihood (ML) setting, it is quite troublesome to perform inference on such quantities, while in
the Bayesian setting it is usually straightforward to obtain the posterior distribution of any non-linear function of the model parameters. Furthermore, in the classical approach, models are usually compared by any other means
than the likelihood. In the Bayesian setting, marginal likelihoods and Bayes factors allow for consistent comparison of non-nested models while incorporating Occam's razor for parsimony. Also, Bayesian estimation provides
reliable results even for finite samples. Finally, \cite{Hall2003} add that the ML approach presents some limitations when the errors are heavy tailed, also the convergence rate is slow and the estimators may not be
asymptotically Gaussian.

This survey reviews the existing Bayesian inference methods for univariate and multivariate GARCH models while having in mind their error specifications. The main emphasis of the paper is on the recent development of an
alternative inference approach for these models using Bayesian non-parametrics. The classical parametric modeling, relying on a finite number of parameters, although so widely used, has some certain drawbacks. Since the number of parameters for any model is fixed, one can encounter underfitting or overfitting, which arises from the misfit between the data available and the parameters needed to estimate. Then, in order to avoid assuming wrong parametric
distributions, which may lead to inconsistent estimators, it is better to consider a semi- or non-parametric approach. Bayesian non-parametrics may lead to less constrained models than classical parametric Bayesian statistics and provide an adequate description of the data, especially when the conditional return distribution is far away from Gaussian.

Up to our knowledge, there has been a very few papers using Bayesian non-parametrics for GARCH models. These are \cite{Ausin2014} for univariate GARCH, \cite{Jensen2013} and \cite{Virbickaite2013} for MGARCH. All of them have considered infinite mixtures of Gaussian distributions with a Dirichlet process (DP) prior over the mixing distribution, which results into DP mixture (DPM) models. This approach so far proves to be the most popular Bayesian
non-parametric modeling procedure. The results over the papers have been consistent: The Bayesian non-parametric approach leads to more flexible models and is better in explaining heavy-tailed return distributions, which
parametric models cannot fully capture.

The outline of this survey is as follows. Section \ref{section:univ_garch} shortly introduces univariate GARCH models and different inference and prediction methods. Section \ref{section:mult_garch} overviews the existing models for multivariate GARCH and different inference and prediction approaches. Section \ref{section:non_param} introduces the Bayesian non-parametric modeling approach and reviews the limited literature of this area in
time-varying volatility models. Section \ref{section:illustration} presents a real data application. Finally, Section \ref{section:conclusions} concludes.

\section{Univariate GARCH}
\label{section:univ_garch}

As mentioned before, the two most popular approaches to model volatility are GARCH type and SV type models. In this survey we focus on GARCH models, therefore, SV models will not be included thereafter. Also, we are not going to enter into the technical details of the Bayesian algorithms and refer to \cite{Robert2004} for a more detailed description of Bayesian techniques.

\subsection{Description of Models}

The general structure of an asset return series modeled by a GARCH-type models can be written as:
\begin{equation*}
r_t=\mu_t+a_t=\mu_t+\sqrt{h_t}\epsilon_t,
\end{equation*}
where $\mu_t=\textrm{E}\left[r_t|\mathcal{I}_{t-1}\right]$ is the conditional mean given ${\mathcal{I}_{t-1}}$, the information up to time ${t-1}$, $a_t$ is the mean corrected returns of the asset at time ${t}$, $h_t=\textrm{Var}\left[r_t|\mathcal{I}_{t-1}\right]$ is the conditional variance given ${\mathcal{I}_{t-1}}$ and ${\epsilon_t}$ is the standard white noise shock. There are several ways to model the conditional mean, ${\mu_t}$.
The usual assumptions are to consider that the mean is either zero, equal to a constant ($\mu_t=\mu$), or follows an ARMA($p$,$q$) process. However, sometimes the mean is also modeled as a function of the variance, say
${g(h_t)}$, which leads to the GARCH-in-Mean models. On the other hand, the conditional variance, ${h_t}$, is usually modeled  using the GARCH-family models. In the basic GARCH model the conditional variance of the returns
depends on a sum of three parts: a constant variance as the long-run average, a linear combination of the past conditional variances and a linear combination of the past mean squared returns. For instance, in the GARCH(1,1)
model, the conditional variance at time $t$ is given by $h_t=\omega+\alpha a^2_{t-1}+\beta h_{t-1}$, for $t=1,\ldots,T$. There are some restrictions which have to be imposed such as ${\omega>0}$, ${\alpha,\beta\geq 0}$ for
positive variance, and ${\alpha+\beta<1}$ for the covariance stationarity.

\cite{Nelson1991} proposed the exponential GARCH (EGARCH) model that acknowledges the existence of asymmetry in the volatility response to the changes in the returns, sometimes also called the "leverage effect", introduced by
\cite{Black1976}. Negative shocks to the returns have a stronger effect on volatility than positive. Other ARCH extensions that try to incorporate the leverage effect are the GJR model by \cite{Glosten1993} and the TGARCH of
\cite{Zakoian1994}, among many others. As \cite{Engle2004} puts it, ``there is now an alphabet soup'' of ARCH family models, such as AARCH, APARCH, FIGARCH, STARCH etc, which try to incorporate such return features as fat tails, volatility clustering and volatility asymmetry. Papers by \cite{Bollerslev1992}, \cite{Bollerslev1994}, \cite{Engle2002}, \cite{Ishida2002} provide extensive reviews of the existing ARCH-type models. \cite{Bera1993} review ARCH type models, discuss their extensions, estimation and testing, also numerous applications. Also, one can find an explicit review with examples and applications concerning GARCH-family models in \cite{Tsay2005} and Chapter 1 in \cite{Terasvirta2009}.

\subsection{Inference Methods}

The main estimation approach for GARCH-family models is the classical maximum likelihood method. However, recently there has been a rapid development of Bayesian estimation techniques, which offer some advantages compared to the frequentist approach as already discussed in the Introduction. In addition, in the empirical finance setting, the frequentist approach presents an uncertainty problem. For instance, optimal allocation is greatly affected by the
parameter uncertainty, which has been recognized in a number of papers, see \cite{Jorion1986} and \cite{Greyserman2006}, among others. These authors conclude that in the frequentist setting the estimated parameter values are
considered to be the true ones, therefore, the optimal portfolio weights tend to inherit this estimation error. However, instead of solving the optimization problem on the basis of the choice of unique parameter values, the
investor can choose the Bayesian approach, because it accounts for parameter uncertainty, as seen in \cite{Kang2011} and \cite{Jacquier2012}, for example. A number of papers in this field have explored different Bayesian
procedures for inference and prediction and different approaches to modeling the fat-tailed errors and/or asymmetric volatility. The recent development of modern Bayesian computational methods, based on Monte Carlo
approximations and MCMC methods have facilitated the usage of Bayesian techniques, see e.g.\ \cite{Robert2004}.

The standard Gibbs sampling procedure does not make the list because it cannot be used due to the recursive nature of the conditional variance: the conditional posterior distributions of the model parameters are not of a simple form. One of the alternatives is the \textit{Griddy-Gibbs} sampler as in \cite{Bauwens1998}. They discuss that previously used importance sampling and Metropolis algorithms have certain drawbacks, such as that they require a careful choice of a good approximation of the posterior density. The authors propose a Griddy-Gibbs sampler which explores analytical properties of the posterior density as much as possible. In this paper the GARCH model has Student-t errors, which allows for fat tails. The authors choose to use flat (uniform) priors on parameters ${(\omega,\alpha,\beta)}$ with whatever region is needed to ensure the positivity of variance, however, the flat prior for the degrees of freedom cannot be used, because then the posterior density is not integrable. Instead, they choose a half-right side of Cauchy. The posteriors of the parameters were found to be skewed, which is a disadvantage for the commonly used Gaussian approximation. On the other hand, \cite{Ausin2007} modeled the errors of a GARCH model with a mixture of two Gaussian distributions. The advantage of this approach compared to that of Student-t errors, is that if the number of the degrees of freedom is very small (less than 5), some moments may not exist. The authors have chosen flat priors for all the parameters, and discovered that there is little sensitivity to the change in the prior distributions (from uniform to Beta), unlike in \cite{Bauwens1998}, where the sensitivity for the prior choice for the degrees of freedom is high. More articles using a Griddy-Gibbs sampling approach are by \cite{Bauwens2002}, who have modeled asymmetric volatility with Gaussian innovations and have used uniform priors for all the parameters, and by \cite{Wago2004}, who explored an asymmetric GARCH model with Student-t errors.

Another MCMC algorithm used in estimating GARCH model parameters, is the \textit{Metropolis-Hastings} (MH) method, which samples from a candidate density and then accepts or rejects the draws depending on a certain acceptance
probability. \cite{Ardia2006} modeled the errors as Gaussian distributed with zero mean and unit variance while the priors are chosen as Gaussian and a MH algorithm is used to draw samples from the joint posterior distribution. The author has carried out a comparative analysis between ML and Bayesian approaches, finding, as in other papers, that some posterior distributions of the parameters were skewed, thus warning against the abusive use of the
Gaussian approximation. Also, \cite{Ardia2006} has performed a sensitivity analysis of the prior means and scale parameters and concluded that the initial priors in this case are vague enough. This approach has been also used by \cite{Muller1998}, \cite{Nakatsuma2000} and \cite{Vrontos2000}, among others. A special case of the MH method is the random walk Metropolis-Hastings (RWMH) where the proposal draws are generated by randomly perturbing the current value using a spherically symmetric distribution. A usual choice is to generate candidate values from a Gaussian distribution where the mean is the previous value of the parameter and the variance can be calibrated to achieve the desired acceptance probability. This procedure is repeated at each MCMC iteration. \cite{Ausin2007} have also carried out a comparison of estimation approaches, Griddy-Gibbs, RWMH and ML. Apparently, RWMH has difficulties in exploring the tails of the posterior distributions and ML estimates may be rather different for those parameters where posterior distributions are skewed.

In order to select one of the algorithms, one might consider some criteria, such as fast convergence for example. \cite{Asai2006} numerically compares some of these approaches in the context of GARCH. The Griddy-Gibbs method is
capable in handling the shape of the posterior by using smaller MCMC outputs comparing with other methods, also, it is flexible regarding parametric specification of a model. However, it can require a lot of computational time.
This author also investigates MH, adaptive rejection Metropolis sampling (ARMS), proposed by \cite{Gilks1995}, and acceptance-rejection MH algorithms (ARMH), proposed by \cite{Tierney1994}. For more in detail about each method in GARCH models see \cite{Nakatsuma2000} and \cite{Kim1998}, among others. Using simulated data, \cite{Asai2006} calculated geometric averages of inefficiency factors for each method. Inefficiency factor is just an inverse of \cite{Geweke1992} efficiency factor. According to this, the ARMH algorithm performed the best. Also, computational time was taken into consideration, where ARMH clearly outperformed MH and ARMS, while Griddy-Gibbs stayed just a bit behind. The author observes that even though the ARMH method showed the best results, the posterior densities for each parameter did not quite explore the tails of the distributions, as desired. In this case Griddy-Gibbs performs better; also, it requires less draws than ARMH. \cite{Bauwens1998} investigate one more convergence criteria, proposed by \cite{Yu1998}, which is based on cumulative sum (cumsum) statistics. It basically shows that if
MCMC is converging, the graph of a certain cumsum statistic against time should approach zero. Their employed Griddy-Gibbs algorithm converged in all four parameters quite fast. Then, the authors explored the advantages and disadvantages of alternative approaches: the importance sampling and the MH algorithm. Considering importance sampling, one of the main disadvantages, as mentioned before, is to find a good approximation of the posterior density (importance function). Also, comparing with Griddy-Gibbs algorithm, the importance sampling requires much more draws to get smooth graphs of the marginal densities. For the MH algorithm, same as in importance sampling, a good approximation needs to be found. Also, compared to Griddy-Gibbs, the MH algorithm did not fully explore  the tails of the distribution, unless for a very big number of draws.

Another important aspect of the Bayesian approach, as commented before, are the advantages in model selection compared to the classical methods. \cite{Miazhynskaia2006} reviews some Bayesian model selection methods using MCMC for GARCH-type models, which allow for the estimation of either marginal model likelihoods, Bayes factors or posterior model probabilities. These are compared to the classical model selection criteria showing that the Bayesian approach clearly considers model complexity in a more unbiased way. Also, \cite{Chen2008} includes a revision of Bayesian selection methods for asymmetric GARCH models, such as the GJR-GARCH and threshold GARCH. They show how
using the Bayesian approach it is possible to compare complex and non-nested models to choose for example between GARCH and stochastic volatility models, between symmetric or asymmetric GARCH models or to determine the number of regimes in threshold processes, among others.

An alternative approach to the previous parametric specifications is the use of Bayesian non-parametric methods, that allow to model the errors as an infinite mixture of normals, as seen in the paper by \cite{Ausin2014}. The Bayesian non-parametric approach for time-varying volatility models will be discussed in detail in Section \ref{section:non_param}.

To sum up, considering the amount of articles published quite recently regarding the topic of estimating univariate GARCH models using MCMC methods indicates still growing interest in the area. Although numerous GARCH-family models have been investigated using different MCMC algorithms, there are still a lot of areas that need further research and development.

\section{Multivariate GARCH}\label{section:mult_garch}

Returns and volatilities depend on each other, so multivariate analysis is a more natural and useful approach. The starting point of multivariate volatility models is a univariate GARCH, thus the most simple MGARCH models can be viewed as direct generalizations of their univariate counterparts. Consider a multivariate return series $\left\{r_{t}\right\}_{t=1}^T$ of size $K\times1$. Then
\begin{equation*}
r_{t} =\mu_{t}+a_{t}=\mu_{t}+H_{t}^{1/2}\epsilon_t,
\end{equation*}
where $\mu_{t}=\tx{E}[r_{t}|\mathcal{I}_{t-1}]$, ${a_t}$ are mean-corrected returns, ${\epsilon_t}$ is a random vector, such that ${\tx{E}[\epsilon_t]=0}$ and ${\tx{Cov}[{\epsilon_t}]=I_K}$ and ${H_t^{1/2}}$ is a positive definite matrix of dimensions ${K\times K}$, such that ${H_t}$ is the conditional covariance matrix of ${r_t}$, i.e., $\tx{Cov}[r_t|\mathcal{I}_{t-1}]=H_{t}^{1/2} \tx{Cov}[\epsilon_t] (H_{t}^{1/2})'=H_t$. There is a wide range of MGARCH models, where most of them differ in specifying ${H_t}$. In the rest of this section we will review the most popular and widely used and the different Bayesian approaches to make inference and prediction. For general
reviews on MGARCH models, see \cite{Bauwens2006a}, \cite{Silvennoinen2008} and \cite{Tsay2005} (Chapter 10), among others.

Regarding inference, one can also consider the same arguments provided in the univariate GARCH case above. Maximum likelihood estimation for MGARCH models can be obtained by using numerical optimization algorithms, such as Fisher scoring and Newton-Raphson. \cite{Vrontos2003} have estimated several bivariate ARCH and GARCH models and found that some classical estimates of the parameters were quite different from their Bayesian counterparts. This was due
to the non-Normality of the parameters. Thus, the authors suggest careful interpretation of the classical estimation approach. Also, \cite{Vrontos2003} found it difficult to evaluate the classical estimates under the stationarity conditions, and consequently the resulting parameters, evaluated ignoring the stationarity constraints, produced non-stationary estimates. These difficulties can be overcome using the Bayesian approach.

\subsection{VEC, DVEC and BEKK}

The VEC model was proposed by \cite{Bollerslev1988}, where every conditional variance and covariance (elements of the ${H_t}$ matrix) is a function of all lagged conditional variances and covariances, as well as lagged squared mean-corrected returns and cross-products of returns. Using this unrestricted VEC formulation, the number of parameters increases dramatically. For example, if ${K=3}$, the number of parameters to estimate will be 78, and if ${K=4}$, the number of parameters increases to 210, see \cite{Bauwens2006a} for the explicit formula for the number of parameters in VEC models. To overcome this difficulty, \cite{Bollerslev1988} simplified the VEC model by proposing a diagonal VEC model, or DVEC, as follows:
\begin{align*}
H_t=\Omega+A\odot (a_{t-1}a_{t-1}')+B\odot H_{t-1},
\end{align*}
where ${\odot}$ indicates the Hadamard product, ${\Omega}$, ${A}$ and ${B}$ are symmetric ${K\times K}$ matrices. As noted in \cite{Bauwens2006a}, ${H_t}$ is positive definite provided that ${\Omega}$, ${A}$, ${B}$ and the
initial matrix ${H_0}$ are positive definite. However, these are quite strong restrictions on the parameters. Also, DVEC model does not allow for dynamic dependence between volatility series. In order to avoid such strong
restrictions on the parameter matrices, \cite{Engle1995} propose the BEKK model, which is just a special case of a VEC and, consequently, less general. It has the attractive property that the conditional covariance matrices are positive definite by construction. The model looks as follows:
\begin{align}\label{bekk}
H_t=\Omega^* \Omega^{*'}+A^* (a_{t-1}a_{t-1}') A^{*'}+B^* H_{t-1} B^{*'},
\end{align}
where ${\Omega^*}$ is a lower triangular matrix and ${A^*}$ and ${B^*}$ are ${K\times K}$ matrices. In the BEKK model it is easy to impose the definite positiveness of the ${H_t}$ matrix. However, the parameter matrices ${A^*}$ and ${B^*}$ do not have direct interpretations since they do not represent directly the size of the impact of the lagged values of volatilities and squared returns.

\cite{Osiewalski2004} present a paper that compares the performance of various bivariate ARCH and GARCH models, such as VEC, BEKK, etc, estimated using Bayesian techniques. As the authors observe, they are the first to perform model comparison using Bayes factors and posterior odds in the MGARCH setting. The algorithm used for parameter estimation and inference is Metropolis-Hastings, and to check for convergence they rely on the cumsum statistics, introduced by \cite{Yu1998}, and used by \cite{Bauwens1998} in the univariate GARCH setting. Using the real data the authors found that the t-BEKK models performed the best, leaving t-VEC not so far behind; t-VEC model, sometimes also called t-VECH, is a more general form of a DVEC, seen above, where the mean-corrected returns follow a Student-t distribution. The name comes from a function called ${vech}$, which reshapes the lower triangular portion of a symmetric variance-covariance matrix into a column vector. To sum up, the authors choose t-BEKK model as clearly better than the t-VEC, because it is relatively simple and has less parameters to estimate.

On the other hand, \cite{Hudson2008} developed a prior distribution for a VECH specification that directly satisfy both necessary and sufficient conditions for positive definiteness and covariance stationarity, while remaining diffuse and non-informative over the allowable parameter space. These authors employed MCMC methods, including Metropolis-Hastings, to help enforce the conditions in this prior.

More recently, \cite{Burda2013b} use the BEKK-GARCH model to show the usefulness of a new posterior sampler called the Adaptive Hamiltonian Monte Carlo (AHMC). Hamiltonian Monte Carlo (HMC) is a procedure to sample from complex distributions. The AHMC is an alternative inferential method based on HMC that is both fast and locally adaptive. The
AHMC appears to work very well when the dimension of the parameter space is very high. Model selection based on marginal likelihood is used to show that full BEKK models are preferred to restricted diagonal specifications. Additionally, \cite{Burda2013} suggests an approach called Constrained Hamiltonian Monte Carlo (CHMC) in order to deal with high dimensional BEKK models with targeting, which allow for a parameter dimension reduction without compromising the model fit, unlike the diagonal BEKK. Model comparison of the full BEKK and the BEKK with targeting is performed indicating that the latter dominates the former in terms of marginal likelihood.

\subsection{Factor-GARCH}

Factor-GARCH was first proposed by \cite{Engle1990} to reduce the dimension of the multivariate model of interest using an accurate approximation of the multivariate volatility. The
definition of the Factor-GARCH model, proposed by \cite{Lin1992}, says that BEKK model in \eqref{bekk} is a Factor-GARCH, if ${A^*}$ and ${B^*}$ have rank one and the same left and right eigenvalues: ${A^*=\alpha w \lambda'}$,
${B^*=\beta w \lambda'}$, where ${\alpha}$ and ${\beta}$ are scalars and ${w}$ and ${\lambda}$ are eigenvectors. Several variants of the factor model have been proposed. One of them is the full-factor multivariate GARCH by \cite{Vrontos2003a}:
\begin{align*}
r_t&=\mu+a_t,\\
a_t&=WX_t,\\
X_t|\mathcal{I}_{t-1}&\sim \mathcal{N}_K (0, \Sigma_t),
\end{align*}
where ${\mu}$ is a ${K\times 1}$ vector of constants, which is time invariant, ${W}$ is a ${K\times K}$ parameter matrix,  ${X_t}$ is a ${K\times 1}$ vector of factors and ${\Sigma_t=\textrm{diag}(\sigma_{1t}^2,\hdots,\sigma_{Kt}^2)}$ is a ${K\times K}$ diagonal variance matrix such that $\sigma_{it}^2=c_i+b_i x_{i,t-1}^2+ g_i \sigma_{i,t-1}^2$, where ${\sigma_{it}^2}$ is the conditional variance
of  the ${i}$th factor at time ${t}$ such that ${c_i>0}$, ${b_i\geq 0}$, ${g_i\geq0}$. Then, the factors in the ${X_t}$ vector are GARCH(1,1) processes and the vector ${a_t}$ is a linear combination of such factors. It can be
easily shown that ${H_t}$ is always positive definite by construction. However, the structure of ${H_t}$ depends on the order of the time series in ${r_t}$. \cite{Vrontos2003a} have considered this problem to find the best
ordering under the proposed model. Furthermore, \cite{Vrontos2003a} investigate a full-factor MGARCH model using the ML and Bayesian approaches. The authors compute maximum likelihood estimates using Fisher scoring algorithm. As for the Bayesian analysis, the authors have adopted a Metropolis-Hastings algorithm, and found that the algorithm is very time consuming, especially in high-dimensional data. To speed-up the convergence, \cite{Vrontos2003a} have proposed reparametrization of positive parameters and also a blocking sampling scheme, where the parameter vector is divided into three blocks: mean, variance and the matrix of constants ${W}$. As mentioned before, the ordering of the univariate time series in full-factor models is important, thus to select ``the best'' model one has to consider ${K!}$ possibilities for a multivariate dataset of dimension ${K}$. Instead of choosing one model and making inference (as if the selected model was the true one), the authors employ a  Bayesian approach by calculating the posterior probabilities for all competing models and model averaging to provide ``combined'' predictions. The main contribution of this paper is that the authors were able to carry out an extensive Bayesian analysis of a full-factor MGARCH model considering not only parameter uncertainty, but model uncertainty as well.

As already discussed above, a very common stylized feature of financial time series is the asymmetric volatility. \cite{Dellaportas2007} have proposed a new class of tree structured MGARCH models that explore the asymmetric
volatility effect. Same as the paper by  \cite{Vrontos2003a}, the authors consider not only parameter-related uncertainty, but also uncertainty corresponding to model selection. Thus in this case Bayesian approach becomes
particularly useful because an alternative method based on maximizing the pseudo-likelihood is only able to work after selecting a single model. The authors develop an MCMC stochastic search algorithm that generates candidate
tree structures and their posterior probabilities. The proposed algorithm converged fast. Such modeling and inference approach leads to more reliable and more informative results concerning model-selection and individual
parameter inference.

There are more models that are nested in BEKK, such as the Orthogonal GARCH for example, see \cite{Alexander1997} and \cite{VanderWeide2002}, among others. All of them fall into the class of direct generalizations of univariate
GARCH or linear combinations of univariate GARCH models. Another class of models are the nonlinear combinations of univariate GARCH models, such as conditional correlation (CCC), dynamic condition correlation (DCC), general
dynamic covariance (GDC) and Copula-GARCH models. A very recent alternative approach that also considers Bayesian estimation can be found in \cite{Jin2013} who proposes a new dynamic component models of returns and realized
covariance (RCOV) matrices based on time-varying Wishart distributions. In particular, Bayesian estimation and model comparison is conducted with a existing range of multivariate GARCH models and RCOV models.

\subsection{CCC}

The CCC model, proposed by \cite{Bollerslev1990} and the simplest in its class, is based on the decomposition of the conditional covariance matrix into conditional standard deviations and correlations. Then, the conditional
covariance matrix ${H_t}$ looks as follows:
\begin{equation*}
H_t=D_t R D_t,
\end{equation*}
where $D_t$ is diagonal matrix with the $K$ conditional standard deviations and ${R}$ is a time-invariant conditional correlation matrix such that ${R=(\rho_{ij})}$ and ${\rho_{ij}=1, \forall i=j}$. The CCC approach can be applied to a wide range of univariate GARCH family models, such as exponential GARCH or GJR-GARCH, for example.

\cite{Vrontos2003} have estimated some real data using a variety of bivariate ARCH and GARCH models in order to select the best model specification and to compare the Bayesian parameter estimates to those of the ML. These authors have considered three ARCH and three GARCH models, all of them with constant conditional correlations (CCC). They have used a Metropolis-Hastings algorithm, which allows to simulate from the joint posterior distribution of the
parameters. For model comparison and selection, \cite{Vrontos2003} have obtained predictive distributions and assessed comparative validity of the analyzed models, according to which the CCC model with diagonal covariance matrix performed the best considering one-step-ahead predictions.

\subsection{DCC}

A natural extension of the simple CCC model are the dynamic conditional correlation (DCC) models, firstly proposed by \cite{Tse2002} and \cite{Engle2002a}. The DCC approach is more realistic, because the dependence between returns is likely to be time-varying.

The models proposed by \cite{Tse2002} and \cite{Engle2002a} consider that the conditional covariance matrix ${H_t}$ looks as $H_t=D_t R_t D_t$, where $R_t$ is now a time varying correlation matrix at time $t$. The models differ in the specification of $R_t$. In the paper by \cite{Tse2002}, the conditional correlation matrix is $R_t=(1-\theta_1-\theta_2)R+\theta_1 R_{t-1}+\theta_2 \Psi_{t-1}$, where ${\theta_1}$ and ${\theta_2}$ are non-negative scalar parameters, such that ${\theta_1+\theta_2< 1}$, ${R}$ is a positive definite matrix such that ${\rho_{ii}=1}$ and ${\Psi_{t-1}}$ is a ${K \times K}$ sample correlation matrix of the past ${m}$ standardized mean-corrected returns ${u_{t}=D_t^{-1}a_t}$.  On the other hand, in the paper by \cite{Engle2002a}, the specification of ${R_t}$ is $R_t=(I\odot Q_t)^{-1/2} Q_t (I\odot Q_t)^{-1/2}$, where $Q_t=(1-\alpha-\beta)\bar Q+\alpha(u_{t-1} u'_{t-1})+\beta Q_{t-1}$, ${u_{i,t}=a_{i,t}/\sqrt{h_{ii,t}}}$ is a mean-corrected standardized returns, ${\alpha}$ and ${\beta}$ are non-negative scalar parameters, such that ${\alpha+\beta<1}$ and ${\bar Q}$ is unconditional covariance matrix
of ${u_t}$. As noted in \cite{Bauwens2006a}, the model by \cite{Engle2002a} does not formulate the conditional correlation as a weighted sum of past correlations, unlike in the DCC model by \cite{Tse2002}, seen above. The
drawback of both these models is that ${\theta_1}$, ${\theta_2}$, ${\alpha}$ and ${\beta}$ are scalar parameters, so all conditional correlations have the same dynamics. However, as \cite{Tsay2005} notes it, the models are parsimonious.

Moreover, as financial returns display not only asymmetric volatility, but also excess kurtosis, previous research, as in univariate case, has mostly considered using a multivariate Student-t distribution for the errors. However, as already discussed above, this approach has several limitations. \cite{Galeano2010} propose a MGARCH-DCC model, where the standardized innovations follow a mixture of Gaussian distributions. This allows to capture long tails  without being limited by the degrees of freedom constraint, which is necessary to impose in the Student-t distribution so that higher moments could exist. The authors estimate the proposed model using the classical ML and Bayesian approaches. In order to estimate model parameters, dynamics of single assets and dynamic correlations, and the parameters of the Gaussian mixture, \cite{Galeano2010} have relied on RWMH algorithm. BIC criteria was used for selecting the number of mixture components, which performed well in simulated data. Using real data, the authors provide an application to calculating the Value at Risk (VaR) and solving a portfolio selection problem. MLE and Bayesian approaches have performed similarly in point estimation, however, the Bayesian approach, besides giving just point estimates, allows the derivation of predictive distributions for the portfolio VaR.

An extension of the DCC model of \cite{Engle2002a} is the Asymmetric DCC also proposed by \cite{Engle2002a}, which incorporates an asymmetric correlation effect. It means that correlations between asset returns decrease more in the bear market than they increase when the market performs well. \cite{Cappiello2006} generalizes the ADCC model into the AGDCC model, where the parameters of the correlation equation are vectors, and not scalars. This allows for asset-specific correlation dynamics. In the AGDCC model, the $Q_t$ matrix in the DCC model is replaced with:
\begin{equation*}
Q_t=S(1-\bar{\kappa}^2-\bar{\lambda}^2-\bar{\delta}^2/2)+\kappa\kappa' \odot u_{t-1}'u_{t-1}+\lambda\lambda'\odot Q_{t-1}+\delta\delta'\odot \eta_{t-1}'\eta_{t-1},
\end{equation*}
where $u_t =D_t^{-1}a_t$ are mean corrected standardized returns, $\eta_t=u_t\odot I(u_t<0)$ selects just negative returns, "diag" stands for either taking just the diagonal elements from the matrix, or making a diagonal matrix from a vector, $S$ is a sample correlation matrix of $u_t$, $\kappa, \lambda$ and $\delta$ are $K\times 1$ vectors, $\bar{\kappa}=K^{-1}\sum_{i=1}^K \kappa_i$, $\bar{\lambda}=K^{-1}\sum_{i=1}^K \lambda_i$ and $\bar{\delta}=K^{-1}\sum_{i=1}^K \delta_i$. To ensure positivity and stationarity of $Q_t$, it is necessary to impose $\kappa_i, \lambda_i, \delta_i > 0$ and $\kappa_i^2+\lambda_i^2+\delta_i^2/2<1$, $\forall i=1,\hdots,K$. The AGDCC by \cite{Cappiello2006} is just a special case where $\kappa_1=\hdots=\kappa_K $, $\lambda_1=\hdots=\lambda_K$ and $\delta_1=\hdots=\delta_K$.

Up to our knowledge, the only paper that considers the AGDCC model in the Bayesian setting is \cite{Virbickaite2013} that propose to model the distribution of the standardized returns as an infinite scale mixture of Gaussian distributions by relying on Bayesian non-parametrics. This approach is presented in more detail in Section \ref{section:non_param}.

\subsection{Copula-GARCH}

The use of copulas is an alternative approach to study return time series and their volatilities. The main convenience of using copulas is that individual marginal densities of the returns can be defined separately from their dependence structure. Then, each marginal time series can be modeled using univariate specification and the dependence between the returns can be modeled by selecting an appropriate copula function. A ${K}$-dimensional copula
${C(u_1,\hdots,u_K)}$, is a multivariate distribution function in the unit hypercube ${[0,1]^K}$, with uniform ${[0,1]}$ marginal distributions. Under certain conditions, the Sklar Theorem affirms that (see, \citealp{Sklar1959}), every joint distribution ${F(x_1,\hdots,x_K)}$, whose marginals are given by ${F_1(x_1),\hdots,F_K(x_K)}$, can be written as $F(x_1,\hdots,x_K)=C(F_1(x_1),\hdots,F_K(x_K))$, where ${C}$ is a copula function of $F$, which is unique  if the marginal distributions are continuous.

The most popular approach to volatility modeling through copulas is called the Copula-GARCH model, where univariate GARCH models are specified for each marginal series and the dependence structure between them is described using a copula function. A very useful feature of copulas, as noted by \cite{Patton2009}, is that the marginal distributions of each random variable do not need to be similar to each other. This is very important in modeling time
series, because each of them might be following different distributions. The choice of copulas can vary from a simple Gaussian copula to more flexible ones, such as Clayton, Gumbel, mixed Gaussian, etc. In the existing literature different parametric and non-parametric specifications can be used for the marginals and copula function $C$. Also, the copula function can be assumed to be constant or time varying, as seen in \cite{Ausin2010}, among others.

The estimation for Copula-GARCH models can be performed in a variety of ways. Maximum likelihood is the obvious choice for fully parametric models. Estimation is generally based on a multi-stage method, where firstly the parameters of the marginal univariate distributions are estimated and then used to condition in estimating the parameters of the copula. Another approach is non- or semi-parametric estimation of the univariate marginal distributions followed by a parametric estimation of the copula parameters. As \cite{Patton2006} has showed, the two-stage maximum likelihood approach lead to consistent, but not efficient estimators.

An alternative is to employ a Bayesian approach, as done by \cite{Ausin2010}. The authors have developed a one-step Bayesian procedure where all parameters are estimated at the same time using the entire likelihood function and, provided the methodology, for obtaining optimal portfolio, calculating VaR and CVaR. \cite{Ausin2010} have used a Gibbs sampler to sample from a joint posterior, where each parameter is updated using a RWMH. In order to reduce
computational cost, the model and copula parameters are updated not one-by-one, but rather by blocks, that consist of highly correlated vectors of model parameters.

\cite{Arakelian2012} have also used Bayesian inference for copula-GARCH models. These authors have proposed a methodology for modelling dynamic dependence structure by allowing copula functions or copula parameters to change across time. The idea is to use a threshold approach so these changes, that are assumed to be unknown, do not evolve in time but occur in distinct points. These authors have also employed a RWMH for parameter estimation together
with a Laplace approximation. The adoption of an MCMC algorithm allows the choice of different copula functions and/or different parameter values between two time thresholds. Bayesian model averaging is considered for predicting dependence measures such as the Kendall's correlation. They conclude that the new model performs well and offers a good insight into the time-varying dependencies between the financial returns.

\cite{Hofmann2010} developed Bayesian inference of a multivariate GARCH model where the dependence is introduced by a D-vine copula on the innovations. A D-vine copula is a special case of vine copulas which are very flexible to construct multivariate copulas because it allows to model dependency between pairs of margins individually. Inference is carried out using a two-step MCMC method closely related with the usual two-step maximum likehood procedure for estimating copula-GARCH models. The authors then focus on estimating the VaR of a portfolio that shows asymmetric dependencies between some pairs of assets and symmetric dependency between others.

All the previously introduced methods rely on parametric assumptions for the distribution of the errors. However, imposing a certain distribution can be rather restrictive and lead to underestimated uncertainty about future volatilities, as seen in \cite{Virbickaite2013}. Therefore, Bayesian non-parametric methods become especially useful, since they do not impose any specific distribution on the standardized returns.

\section{Bayesian Non-Parametrics for GARCH}\label{section:non_param}

Bayesian non-parametrics is an alternative approach to the classical parametric Bayesian statistics, where one usually gives some prior for the parameters of interest, whose distribution is unknown, and then observes the data and calculates the posterior. The priors come from the family of parametric distributions. Bayesian non-parametrics uses a prior over distributions with the support being the space of all distributions. Then, it can be viewed as a
distribution over distributions.

\subsection{DP and DPM}

One of the most popular Bayesian non-parametric modeling approach is based on Dirichlet processes (DP) and mixtures of Dirichlet processes (DPM). DP was firstly introduced by  \cite{Ferguson1973}. Suppose that we have a sequence of exchangeably distributed random variables ${X_1, X_2,\hdots}$ from an unknown distribution ${F}$, where the support for ${X_i}$ is ${\Omega}$. In order to perform Bayesian inference, we need to define the prior for ${F}$. This can be done by considering partitions of ${\Omega}$, such as ${\Omega=C_1\cup C_2 \cup \hdots \cup C_m}$, and defining priors over all possible partitions. We say that ${F}$ has a Dirichlet process prior, denoted as ${F\sim\mathcal{DP}(\alpha,F_0)}$, if the set of associated probabilities given ${F}$ for any partition follows a Dirichlet distribution, $
\{F(C_1),\hdots,F(C_m)\}\sim \textrm{Dirichlet}  \:\:(\alpha F_0(C_1),\hdots,\alpha F_0 (C_m))$,
where ${\alpha >0}$ is a precision parameter that represents our prior certainty of how concentrated the distribution is around ${F_0}$, which is a known base distribution on ${\Omega}$. The Dirichlet process is a conjugate prior. Thus, given ${n}$ independent and identically distributed samples from ${F}$, the posterior distribution of ${F}$ is also a Dirichlet process such that $
F\sim\mathcal{DP} \left({\alpha +n},\left(\alpha F_0 +n F_n\right)\left(\alpha +n\right)^{-1}\right),
$
where ${F_n}$ is the empirical distribution function.

There are two main ways for generating a sample from the marginal distribution of ${X}$, where ${X|F\sim F}$ and ${F\sim \mathcal{DP}(\alpha,F_0)}$: the Polya urn and stick breaking procedures. On the one hand, the Polya urn scheme can be illustrated in terms of a urn with ${\alpha}$ black balls; when a non-black ball is drawn, it is placed back in the urn together with another ball of the same color. If the drawn ball is black, a new color is generated from $F_0$ and a ball of this new color is added to the urn together with the black ball we drew. This process gives a discrete marginal distribution for ${X}$ since there is always a probability that a previously seen value is repeated. On the other hand, the stick-breaking procedure is based on the representation of the random distribution ${F}$ as a countably infinite mixture:
\begin{align*}
F=\sum_{m=1}^{\infty} \omega_m \delta_{X_m},
\end{align*}
where ${\delta_X}$ is a Dirac measure, ${X_m\sim F_0}$ and the weights are such that $ \omega_1=\beta_1,\:\:\omega_m=\beta_m \prod_{i=1}^{m-1}(1-\beta_i),\:\:\textrm{for}\:\: m=1,\hdots,$  where ${\beta_m\sim  \textrm{Beta}\:(1,\alpha)}$. This implies that the weights ${\omega \rightarrow \textrm{Dirichlet} (\alpha/K,\hdots,\alpha/K)}$ as ${K \rightarrow \infty}$.

The discreteness of the Dirichlet process is clearly a disadvantage in practice. A solution was proposed by \cite{Antoniak1974} by using DPM models where a DP prior is imposed over over the distribution of the model parameters, $\theta$, as follows:
\begin{align*}
  X_i|\theta_i&\sim F(X|\theta_i),\\
  \theta_i|G&\sim G(\theta),\\
  G|\alpha,G_0&\sim \mathcal{DP} (\alpha,G_0).
\end{align*} Observe that ${G}$ is a random distribution drawn from the DP and because it is discrete, multiple ${\theta_i}$'s can take the same value simultaneously, making it a mixture model. In fact, using the stick-breaking representation, the hierarchical model above can be seen as an infinite mixture of distributions:
\begin{align*}
f(X|\theta,\omega)=\sum\limits_{m=1}^{\infty}\omega_m f(X|\theta_m),
\end{align*}
where the weights are obtained as before: ${\omega_1=\beta_1}$, ${\omega_m=\beta_m \prod_{i=1}^{m-1}(1-\beta_i)}$, for ${m=1,\hdots}$, and where ${\beta_m\sim \textrm{Beta}\:(1,\alpha)}$ and ${\theta_m\sim G_0}$.

Regarding inference algorithms, there are two main types of approaches. On the one hand, the marginal methods, such as those proposed by \cite{Escobar1995}, \cite{Maceachern1994} and \cite{Neal2000}, which rely on the Polya
urn representation. All these algorithms are based on integrating out the infinite dimensional part of the model. Recently, another class of algorithms, called conditional methods, have been proposed. These approaches, based on
the stick-breaking scheme, leave the infinite part in the model and sample a finite number of variables. These include the procedure by \cite{Walker2007}, who introduces slice sampling schemes to deal with the infiniteness in
DPM, and the retrospective MCMC method of \cite{Papaspiliopoulos2008}, that is later combined by \cite{Papaspiliopoulos} with slice sampling method by \cite{Walker2007} to obtain a new composite algorithm, which is better, faster and easier to implement. Generally, the stick-breaking compared to the Polya urn procedures produce better mixing and simpler algorithms.

\subsection{Volatility modeling using DPM}

As mentioned above, so far there has been little research in modeling volatilities with MGARCH using the DPM models. Up to our knowledge, these only include: \cite{Ausin2014} for univariate GARCH, and \cite{Jensen2013} and \cite{Virbickaite2013}, for MGARCH.

\cite{Ausin2014} have applied semi-parametric Bayesian techniques to estimate univariate GARCH-type models. These authors have used the class of scale mixtures of Gaussian distributions, that allow for the variances to change
over components, with a Dirichlet process prior on the mixing distribution to model innovations of the GARCH process. The resulting class of models is called DPM-GARCH type models. In order to perform Bayesian inference on the
new model, the authors employ a stick-breaking sampling scheme and make use of the ideas proposed in \cite{Walker2007}, \cite{Papaspiliopoulos2008} and \cite{Papaspiliopoulos}. The new scale mixture model was compared to a
simpler mixture of two Gaussians, Student-t and the usual Gaussian models. The estimation results in all three cases were quite similar, however, the scale mixture model is able to capture skewness as well as kurtosis and, based on the approximated Log Marginal Likelihood (LML) and $\textrm{DIC}$, provided the best performance in simulated and real data. Finally, \cite{Ausin2014} have applied the resulting model to perform one-step-ahead predictions
for volatilities and VaR. In general, the non-parametric approach leads to wider Bayesian credible intervals and can better describe long tails.

\cite{Jensen2013} propose a Bayesian non-parametric modeling approach for the innovations in MGARCH models. They use a MGARCH specification, proposed by \cite{Ding2001}, which is a different representation of a well known DVEC model, introduced above. The innovations are modeled as an infinite mixture of multivariate Normals with a DP prior. The authors have employed Polya urn and stick-breaking schemes and, using two data sets, compared the three model specifications: parametric MGARCH with Student-t innovations (MGARCH-t), GARCH-DPM-$\Lambda$ that allows for different covariances (scale mixture) and MGARCH-DPM, allowing for different means and covariances of each component (location-scale mixture). In general, both semi-parametric models produced wider density intervals. However, in MGARCH-t model a single degree of freedom parameter determines the tail thickness in all directions of the density, meanwhile the non-parametric models are able to capture various  deviations from Normality by using a certain number of components. These results are consistent with the ones in \cite{Ausin2014}. As for predictions, both semi-parametric models performed equally good and outperformed the parametric MGARCH-t specification.

Finally, the paper by \cite{Virbickaite2013} can be seen as a direct generalization of the paper by \cite{Ausin2014} to the multivariate framework. Here, same as in \cite{Jensen2013}, the authors have proposed using an infinite scale mixture of Normals for the standardized returns. For the MGARCH model a GJR-ADCC was chosen, allowing for asymmetric volatilities and asymmetric time-varying correlations.  Moreover, the authors have carried-out a simulation study that illustrated the adaptability of the DPM model. Finally, the authors provided one real data application to portfolio decision problem concluding that DPM models are less restrictive and more adaptive to whatever distribution the data comes from, therefore, can better capture the uncertainty about financial decisions.

To sum up, the findings in the above papers are consistent: the Bayesian semi-parametric approach leads to more flexible models and is better in explaining heavy-tailed return distributions, which parametric models cannot fully capture. The parameters are less precise, i.e. wider Bayesian credible intervals are observed because the semi-parametric models are less restricted. This provides a more adequate measure of uncertainty. If in the Gaussian
setting the credible intervals are very narrow and the real data is not Gaussian, this makes the agent overconfident about her decisions, and she takes more risk than she would like to assume. \cite{Steel} observes that the
combination of Bayesian methods and MCMC computational algorithms provide new modeling possibilities and calls for more research regarding non-parametric Bayesian time series modeling.

\section{Illustration}\label{section:illustration}

This illustration study using real data has basically two goals: firstly, to show the advantages of the Bayesian approach, such as the ability to obtain posterior densities of quantities of interest and the facility to incorporate various constraints on the parameters. Secondly, to illustrate the flexibility of the Bayesian non-parametric approach for GARCH modeling.

\begin{figure}[tbp]
\centering
\caption{Log-Returns and Histograms of FTSE100 and S\&P500 Indices}
\includegraphics[scale=0.45]{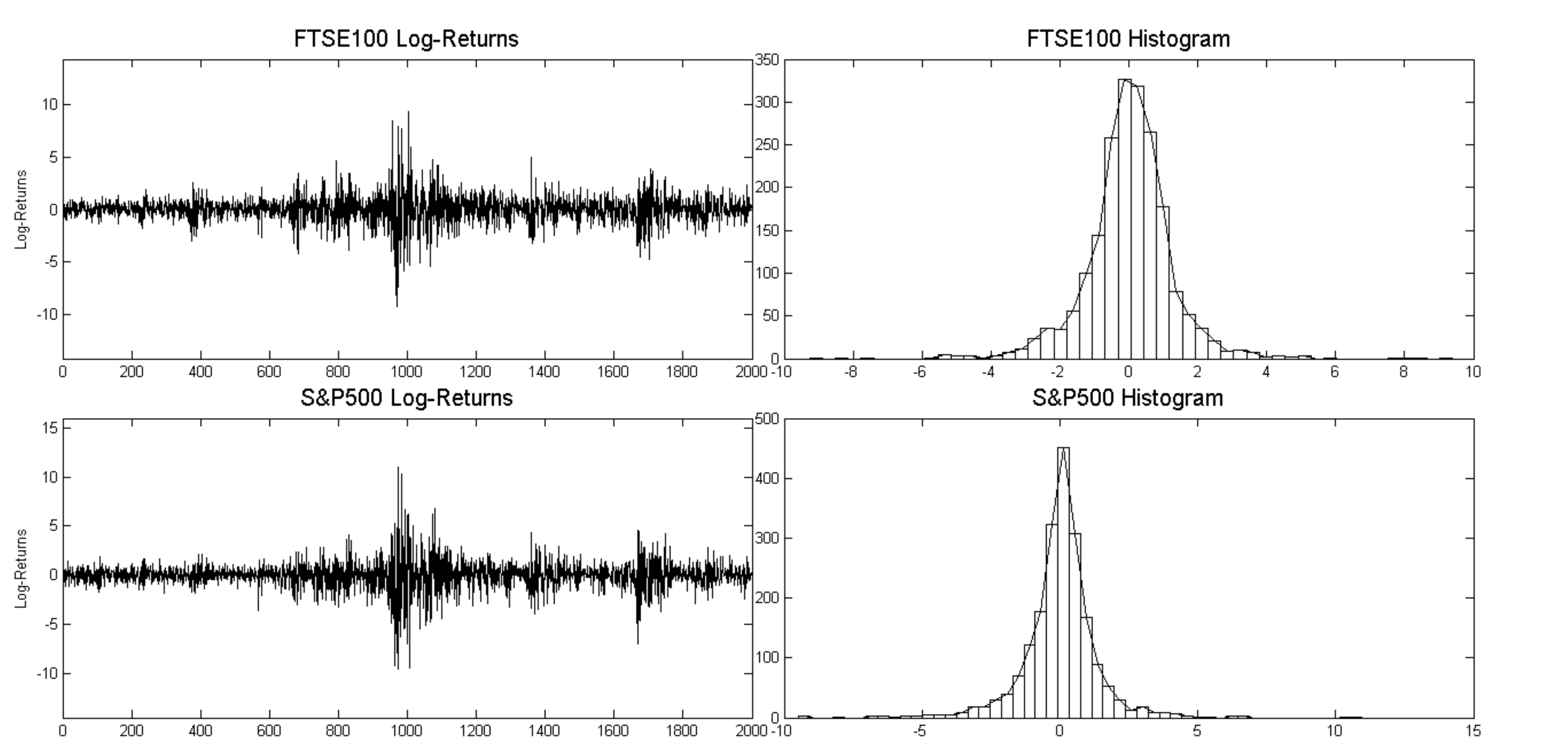}
\label{F:returns}
\end{figure}

The data used for estimation are the log-returns (in percentages), obtained from close prices adjusted for dividends and splits, of two market indices: FTSE100 and S\&P500 from November 10th, 2004 till December 10th, 2012, resulting into a sample size of 2000 observations. FTSE100  is a share index of the 100 companies listed on the London Stock Exchange with the highest market capitalization. S\&P500  is a stock market index based on the common stock prices of 500 top publicly traded American companies. The data was obtained from Yahoo Finance. Figure \ref{F:returns} and Table \ref{T:descriptive_table} present the basic plots and descriptive statistics of the two log-return series.

\begin{table}[tbp]
\centering
\caption{Descriptive Statistics of the FTSE100 and S\&P500 Log-Return Series}
\label{T:descriptive_table}
\begin{small}
\begin{tabular}{ccc}
\hline
\hline\toprule
\rule{0pt}{4ex} & FTSE100 & S\&P500 \\
\hline
Mean &0.0112  &  0.0099\\
Median &0.0344   & 0.0779\\
Variance &1.7164  &  1.9617\\
Skewness &-0.0974  & -0.3001\\
Kurtosis &10.5464  & 12.5674\\
Correlation &0.6060 \\
\hline
\end{tabular}
\end{small}
\end{table}

As seen from the plot and descriptive statistics, the data is slightly skewed and with high kurtosis, therefore, assuming a Gaussian distribution for the standardized returns would be inappropriate. Therefore, we estimate this
bivariate time series using an ADCC model by \cite{Engle2002a}, presented in Section 3.4., which incorporates an asymmetric correlation effect. The univariate series are assumed to follow GJR-GARCH$(1,1)$ models in order to incorporate the leverage effect in volatilities. As for the errors, we use an infinite scale mixture of Gaussian distributions. Therefore, we call the final model GJR-ADCC-DPM. Inference and prediction is carried out using Bayesian non-parametric techniques, as seen in  \cite{Virbickaite2013}. The selection of the  MGARCH specification is arbitrary and other models might work equally well. For the sake of comparison, we estimate a restricted GJR-ADCC-Gaussian model using Maximum Likelihood and Bayesian approaches. The estimation results are presented in Table \ref{table:estimation}.

\begin{table}[tbp]
\caption{Estimation Results for FTSE100 (subindex 1) and S\&P500 (subindex 2) log-Returns, with 30,000 iterations plus 10,000 burn-in iterations}
\centering
\label{table:estimation}
\begin{small}
\begin{tabular}{c c c c c c c}
\hline\hline
  \toprule
  & \multicolumn{2}{c}{ML Gaussian}& \multicolumn{2}{c}{Bayesian Gaussian} & \multicolumn{2}{c}{Bayesian DPM}  \\

&Estimate & St. Dev. & Mean & $95\%$ CI&Mean & $95\%$ CI\\
\\\hline
$\omega_1$ 	& 0.0166 &  0.0020  & 0.0192 & (0.0130, 0.0258) & 0.0181 & (0.0104, 0.0264) \\
$\omega_2$ 	& 0.0190 & 0.0016  & 0.0249 & (0.0174, 0.0316)  & 0.0219 & (0.0153, 0.0293) \\
$\alpha_1$		& $8.31\cdot 10^{-8}$  &$ 2.61\cdot 10^{-4}$  & 0.0058 & (0.0002,    0.0177) &  0.0046  & (0.0003, 0.0112) \\
$\alpha_2$		& $9.05\cdot 10^{-9}$ & $7.58\cdot 10^{-5}$  & 0.0053 & (0.0002,    0.0173) &  0.0059  & (0.0002, 0.0151) \\
$\beta_1$		& 0.9087 &  0.0045  & 0.9010 & (0.8841, 0.9152) & 0.8956 & (0.8762, 0.9139) \\
$\beta_2$		& 0.9079 & 0.0050 & 0.8888 & (0.8705, 0.9088) & 0.8851 & (0.8675, 0.9041) \\
$\phi_1$			& 0.1535 & 0.0085  & 0.1587 & (0.1351, 0.1871) & 0.1586 & (0.1057,  0.2089) \\
$\phi_2$			& 0.1483 &  0.0092  & 0.1737 & (0.1398, 0.2020) & 0.1758 & (0.1134, 0.2142) \\
$\kappa$			& 0.0075 &  0.0020  & 0.0071 & (0.0014, 0.0145) & 0.0095 & (0.0040, 0.0156) \\
$\lambda$		& 0.9898 &  0.0029  & 0.9818 & (0.9665, 0.9898) & 0.9806 & (0.9693, 0.9901) \\
$\delta$			&$5.50\cdot 10^{-8}$  &   $1.33\cdot 10^{-4}$  & 0.0076 & (0.0002,    0.0153) &   0.0039 & (0.0001, 0.0114) \\
\hline
\end{tabular}
\end{small}
\end{table}

The estimated parameters are very similar for all three approaches, except for $\alpha$, the asymmetric correlation coefficient $\delta$. Since $\alpha$ and $\delta$ are so close to zero, the ML has some trouble in estimating those parameters. Overall, the $\delta$ is small, indicating little evidence of asymmetrical behavior in correlations.

\begin{figure}[tbp]
\centering
\caption{Log-Predictive Densities of the one-step-ahead Returns $r_{t+1}$ for Bayesian Gaussian and DPM Models}
\label{fig: tails}
\includegraphics[width=1\textwidth, height=.3\textheight]{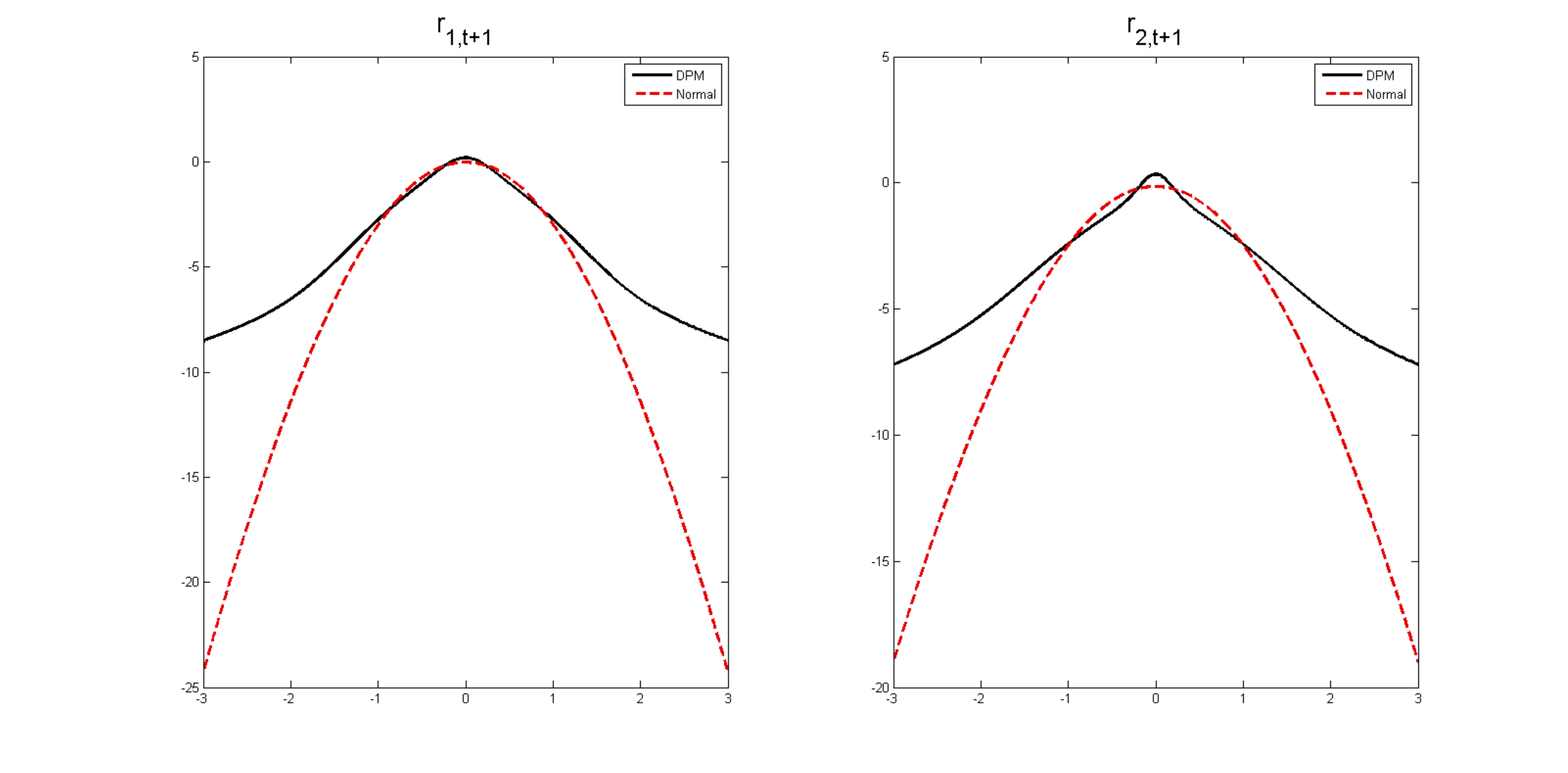}
\end{figure}

Figure \ref{fig: tails} shows the estimated marginal predictive densities of the one-step-ahead returns in log scale using the Bayesian approach. We can observe the differences in tails arising from different specification of the errors. The DPM model allows for a more flexible distribution, therefore, for more extreme returns, i.e.\ fatter tails. The estimated densities were obtained using the procedure described in \cite{Virbickaite2013}.

\begin{figure}[tbp]
\centering
\caption{Densities of One-Step-Ahead Volatilities of the Returns}
\label{fig: vols_returns}
\includegraphics[width=1\textwidth, height=.3\textheight]{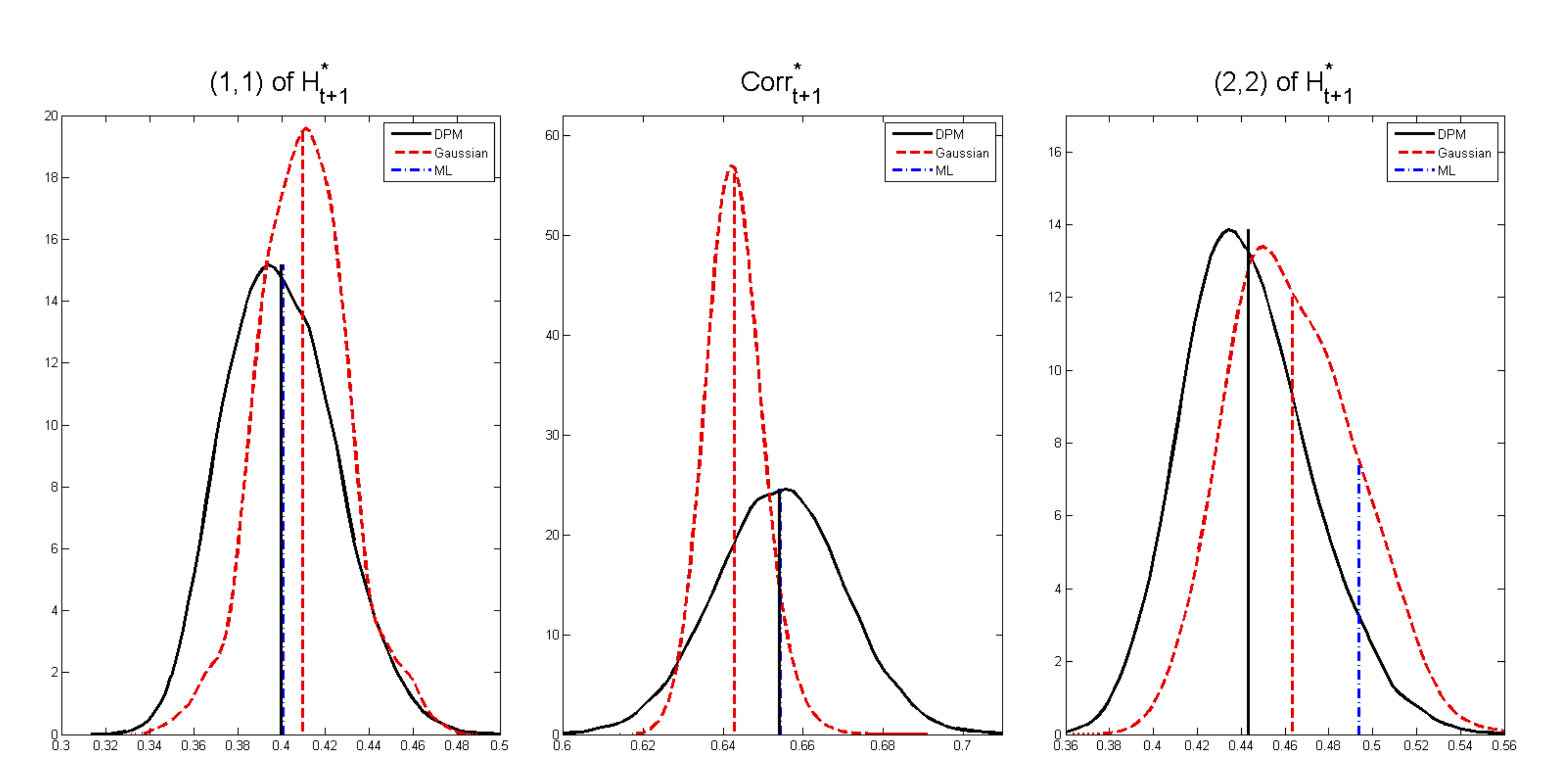}
\end{figure}

Table \ref{table:volatilities} presents the estimated mean, median and $95\%$ credible intervals of one-step-ahead volatility matrices in Bayesian context.  The matrix element (1,1) represents the volatility for the  FTSE100 series, (2,2) for the S\&P500, and the elements in the diagonal (1,2) and (2,1) represent the covariance of both financial returns. Figure \ref{fig: vols_returns} draws the posterior distributions for volatilities and correlation. The estimated mean volatilities for both, DPM and Gaussian approaches, are very similar, however, the main differences arise from the shape of the posterior distribution. $95\%$ credible intervals for DPM model correlation are wider providing a more realistic measure of uncertainty about future correlations between two assets. This is a very important implication in financial setting, because if an investor chooses to be Gaussian, she would be overconfident about her decision and unable to adequately measure the risk she is facing. See \cite{Virbickaite2013} for a more detailed comparison of DPM and alternative parametric approaches in portfolio decision problems.

\begin{table}[tbp]
\caption{Estimated Means, Medians and $95\%$ Credible Intervals of One-Step-Ahead Volatilities of FTSE100 and S\&P500 log-Returns}
\centering
\label{table:volatilities}
\begin{small}
\begin{tabular}{c c c c c c c}
\hline\hline
  \toprule
  & \multicolumn{2}{c}{ }& \multicolumn{2}{c}{Bayesian Gaussian} & \multicolumn{2}{c}{Bayesian DPM}  \\
&Constant &ML Gaussian &Mean & $95\%$ CI&Mean & $95\%$ CI\\
& &&Median & CI Length&Median & CI Length\\
\\\hline
$H_{T+1}^{\star(1,1)}$ & 1.7164 & 0.4007 & 0.4098  & (0.3681, 0.4538) & 0.3996 & (0.3550, 0.4512) \\
&&&  0.4099 & 0.0857 &0.3983& 0.0962 \\
$H_{T+1}^{\star(1,2)}$ & 1.1120 & 0.2911 & 0.2800 & (0.2571, 0.3077 ) & 0.2751 & (0.2421, 0.3123) \\
&&&  0.2790 &0.0506  &0.2742&  0.0702\\
$H_{T+1}^{\star(2,2)}$ & 1.9617 & 0.4939 & 0.4635 & (0.4159, 0.5193) & 0.4431 & (0.3912, 0.5059 ) \\
&&& 0.4606 & 0.1034 &0.4408&  0.1146\\

\hline
\end{tabular}
\end{small}
\end{table}

To sum up, this illustration has shown the main differences between the standard estimation procedures and the new non-parametric approach. Even though the point estimates for the parameters and the one-step-ahead volatilities are very similar, the main differences arise from the thickness of tails of predictive distributions of one-step-ahead returns and the shape of the posterior distribution for the one-step-ahead volatilities.

\section{Conclusions}\label{section:conclusions}

In this paper we reviewed univariate and multivariate GARCH models and inference methods, putting emphasis on the Bayesian approach. We have surveyed the existing literature that concerns various Bayesian inference methods for
MGARCH models, outlining the advantages of the Bayesian approach versus the classical procedures. We have also discussed in more detail the recent Bayesian non-parametric method for GARCH models, which avoid imposing arbitrary parametric distributional assumptions. This new approach is more flexible and can describe better the uncertainty about future volatilities and returns, as has been illustrated using real data.

\section*{Acknowledgements}

We are grateful to an anonymous referee for helpful comments. The first and second authors are grateful for the financial support from MEC grant ECO2011-25706. The third author acknowledges financial support from MEC grant ECO2012-38442.

\bibliographystyle{model5-names}

\end{document}